\documentclass[12pt]{amsart}
\usepackage{amssymb}
\usepackage{amsmath}
\usepackage{pictex}
   \newcommand\al{\alpha}
   \newcommand\AND{\quad\mbox{and}\quad}
   \newcommand\BB{\mathcal B}
   \newcommand\bd{\partial}
   \newcommand\cf{\curlywedge}
   \newcommand\de{\delta}
   \newcommand\df{\operatorname{\sf def}}
   \newcommand\DL{\mbox{\sl DL}}
   \newcommand\dps{\displaystyle}
   \newcommand\Ga{\Gamma}
   \newcommand\geo[1]{\overline{#1}}
   \newcommand\HH{\mathcal H}
   \newcommand\hor{\mathfrak{h}}
   \newcommand\la{\lambda}
   \newcommand\MM{\mathcal M}
   \newcommand\om{\omega}

   \newcommand\R{\mathbb R}
   \newcommand\sib{\,{\buildrel s \over \sim}\,}
   \newcommand\T{\mathbb T}
   \newcommand\wh{\widehat}
   
   \newcommand\Z{\mathbb Z}

\numberwithin{equation}{section}

\newtheoremstyle{mythm}
  {9pt}
  {9pt}
  {\itshape}
  {0pt}
  {\bfseries}
  {}
  { }
  {\thmnumber{(#2)}\thmname{ #1}\thmnote{ #3}}

\newtheoremstyle{mydef}
  {9pt}
  {9pt}
  {\normalfont}
  {0pt}
  {\bfseries}
  {}
  { }
  {\thmnumber{(#2)}\thmname{ #1}\thmnote{ #3}}

\theoremstyle{mythm}
\newtheorem{thm}[equation]{Theorem.}
\newtheorem{pro}[equation]{Proposition.}
\newtheorem{lem}[equation]{Lemma.}
\newtheorem{cor}[equation]{Corollary.}

\theoremstyle{mydef}
\newtheorem{dfn}[equation]{Definition.}
\newtheorem{exa}[equation]{Example.}
\newtheorem{imp}[equation]{}

\begin{document}
\markright{Random walks and harmonic functions}
\title{\large Lamplighters, Diestel-Leader graphs,\\ 
random walks, and harmonic functions}
\author{\bf Wolfgang WOESS}
\email{woess@TUGraz.at}
\address{\parbox{1.4\linewidth}{Institut f\"ur Mathematik C, 
Technische Universit\"at Graz\\
Steyrergasse 30, A-8010 Graz, Austria}}
\date{Version of November 27, 2003, to appear in Combinatorics, 
Probability \& Computing}
\thanks{Partially supported by FWF (Austrian Science Fund) project P15577}
\begin{abstract}
The lamplighter group over $\Z$ is the wreath product $\Z_q \wr \Z$. 
With respect to a natural generating set, its Cayley graph is the 
Diestel-Leader graph $\DL(q,q)$. We study  harmonic functions for the 
``simple'' Laplacian on this graph, and more generally, for a class 
of random walks on $\DL(q,r)$, where $q,r \ge 2$. The $\DL$-graphs 
are horocyclic products of two trees, and we give a full description
of all positive harmonic functions in terms of 
the boundaries of these two trees. In particular, we determine the 
minimal Martin boundary, that is, the set of  minimal positive harmonic 
functions.
\end{abstract}
\maketitle

\markboth{\sf W. Woess}{\sf Random walks and harmonic functions}
\baselineskip 15pt

\section{Introduction}
Think of a (typically infinite) connected graph $X$ where in each vertex 
there is a lamp that may be switched off (state $0$), or switched on with 
$q-1$ different intensities (states $1, \dots, q-1$). 
Initially, all lamps are turned off, and a lamplighter starts at some vertex 
of $X$ and walks around. When he visits a
vertex, he may switch the lamp sitting there into one of its $q$ different states (including ``off'').
Our information consists of the position $x\in X$ of the lamplighter and of the finitely supported 
configuration $\eta: X \to \Z_q = \{0, \dots q-1 \}$ of the lamps that are switched on, 
including their respective intensities. The set $\Z_q \wr X$ of all such pairs $(\eta,x)$ can be equipped
in several ways with a natural connected graph structure, giving rise to a \emph{lamplighter graph}.

When $X$ is a Cayley graph of a group $\Ga$ then underlying this construction, there
is the \emph{wreath product} $\Z_q \wr \Ga$, which is the semidirect product of $\Ga$
with the group of all finitely supported functions $\eta: \Ga \to \Z_q$ (i.e., a direct sum),
on which $\Ga$ acts by $g\eta(h) = \eta(g^{-1}h)$. Instead of $\Z_q = \Z/(q\Z)$,
one may of course take any other group $L$ of ``lamps'', leading to the wreath 
product $L \wr \Ga$.

Various aspects of random walks on lamplighter groups have received considerable attention
recently: Poisson boundary ({\sc Kaimanovich and Vershik} \cite{KaiVer} and 
{\sc Kaimanovich} \cite{Kai}), rate of escape ({\sc Lyons, Pemantle and Peres} 
\cite{LyoPemPer}, {\sc Erschler} \cite{Ers}, {\sc Revelle} \cite{Rev1}), 
spectral theory ({\sc Grigorchuk and \.Zuk} \cite{GriZuk}),
and the asymptotic behaviour of transition probabilites ({\sc Saloff-Coste and
Pittet} \cite{PitSal1}, \cite{PitSal2}, {\sc Revelle} \cite{Rev2}).
Here, we shall consider \emph{harmonic functions.}

A harmonic function on a locally finite graph is a real-valued function 
whose value at each vertex coincides with the arithmetic average of its 
values in the neighbour vertices. More generally, we can consider the 
transition matrix $P$ of a random walk on the graph, suitably adapted to
the graph's geometry; a harmonic function $h$ is then one that satisfies $Ph=h$.

In the present paper, we shall determine all positive harmonic functions 
on certain Cayley graphs of the simplest lamplighter group,  
$\Ga = \Z_q \wr \Z$. Namely, we first explain that the Diestel-Leader 
graph $\DL(q,q)$ is a Cayley graph of $\Ga$. More generally, if $q, r \ge 2$ 
then $\DL(q,r)$ is obtained as a ``horocyclic product'' of two homogeneous 
trees $\T_q$ and $\T_r$  with degrees $q+1$ and $r+1$, respectively. 
We remark that this does not mean that $\DL(q,r)$ is ``almost'' a tree in 
any sense; indeed, it is  a one-ended, vertex-transitive graph
which is a Cayley graph only when $r=q$. When $r \ne q$, it is believed to be 
an example of a transitive graph that is not quasi-isometric with any 
Cayley graph of some finitely generated group -- see {\sc Diestel and Leader} 
\cite{DieLea}.

Nevertheless, we can use the boundary of each of the two trees that 
compose $\DL(q,r)$ for giving an integral representation of all positive 
harmonic functions: in that representation, we start with the projections 
of the random walk on $\DL(q,r)$ to each of the two trees
and the corresponding Martin kernels. Our main result is that every
positive harmonic function on $\DL(q,r)$ is of the form $h = h_1 + h_2$, where 
$h_1$ is obtained by lifting a harmonic function from $\T_q$ to $\DL(q,r)$, 
and $h_2$ is obtained analogously from $\T_r$. 
Thereby, we also determine all minimal positive harmonic functions.

\vspace{.4cm}

We now give an outline of the contents of this paper.

\vspace{.2cm}

Section \ref{geometry}, although it does not contain \emph{proofs,} is
crucial, since it \emph{explains} the geometry of the structures that 
we are working with, and in particular, the correspondence between 
lamplighter groups and Diestel-Leader graphs. As a matter of fact, it is
precisely this geometric realization that allows us
to determine all positive harmonic functions on $\Z_q \wr \Z$.
At the end of \S \ref{geometry}, we state the first main result, regarding
the decomposition of positive harmonic functions over the two trees 
(Theorem \ref{split-theorem}).

In Section \ref{basics}, we recall basic results on positive harmonic 
functions for irreducible Markov chains. In particular, we consider
finite sets with boundaries, the Martin boundary at infinity and its
minimal part, and the Martin compactification for nearest neighbour
random walks on trees. 

In Section \ref{principal}, we use all the preceding ingredients to
prove the Decomposition Theorem \ref{split-theorem}.
It is then quite simple to determine all minimal positive harmonic
functions (Theorem \ref{minimal-theorem}); they are the Martin kernels
of the two projected random walks, up to one, resp. two exceptions.
We then retranslate these results to the lamplighter group $\Z_q \wr \Z$
(Example \ref{SRW-example}).

In Section \ref{extension}, we adapt the preceding results to the 
``switch-walk-switch'' random walk, which is in some sense more natural from the
point of view of the lamplighter than the simple random walk on $\DL(q,q)$.

Section \ref{final} is devoted to some additional remarks and speculations.

\section{Diestel-Leader graphs and lamplighters}\label{geometry}

Let $\T = \T_q$ be the homogeneous tree with degree $q+1$, $q \ge 2$.
A \emph{geodesic path}, resp. \emph{geodesic ray}, resp. \emph{infinite
geodesic} in $\T$ is a finite, 
resp. one-sided infinite, resp. doubly infinite sequence $(x_n)$ of vertices
of $\T$ such that $d(x_i,x_j) = |i-j|$ for all $i, j$, 
where $d(\cdot,\cdot)$ denotes the graph distance. 

Two rays are \emph{equivalent} if their symmetric difference is finite.
An \emph{end} of $\T$ is an equivalence class of rays. The space of 
ends is denoted $\bd \T$, and we write $\wh \T = \T \cup \bd \T$. 
For all $w, z \in \wh \T$ there is a unique geodesic $\geo{w\,z}$ 
that connects the two. In particular, if $x \in \T$ and $\xi \in \bd \T$ then 
$\geo{x\,\xi}$ is the ray that starts at $x$ and represents $\xi$.
Furthermore, if $\xi, \zeta \in \bd \T$ ($\xi \ne \zeta$) then
$\geo{\zeta\,\xi}$ is the infinite geodesic whose two halves (split at any
vertex) are rays that respresent $\zeta$ and $\xi$, respectively.

For $x,y \in \T$, $x \ne y$, we define the cone
$\wh \T(x,y) = \{ w \in \wh \T : y \in \geo{x\,w} \}$.
The collection of all cones is the basis of  a topology wich 
makes $\wh \T$ a compact, totally disconnected Hausdorff space
with $\T$ as a dense, discrete subset. We denote 
$\T(x,y) = \T \cap\wh \T(x,y)$ and $\bd \T(x,y) = \bd \T \cap \T(x,y)$.

We fix a root  $o \in \T$.
If $w, z \in \wh \T$, then their \emph{confluent} $c=w \wedge z$ with
respect to the root vertex $o$ is defined by
$\geo{o\,w} \cap \geo{o\,z} = \geo{o\,c}$. 
Similarly, we choose and fix a reference end $\om \in \bd \T$. For 
$z, v \in \wh \T \setminus \{ \om \}$, their confluent $b = v \cf z$ 
with respect to $\om$ is defined by
$\geo{v\,\om} \cap \geo{z\,\om} = \geo{b\,\om}$.
The \emph{Busemann function} $\hor: \T \to \Z$ and the \emph{horocycles} $H_k$
with respect to $\om$ are defined as
$$
\hor(x) = d(x,x \cf o) - d(o,x \cf o) \AND H_k = \{ x \in \T : \hor(x) = k \}\,.
$$
Every horocycle is infinite. Every vertex $x$ in $H_k$ has one neighbour 
$x^-$ (its predecessor) in $H_{k-1}$ and $q$ neighbours (its successors)
in $H_{k+1}$. We set $\bd^* \T = \bd \T \setminus \{\om\}$. 

\vspace{-.4cm}

$$
\beginpicture 

\setcoordinatesystem units <.7mm,1.04mm>

\setplotarea x from -10 to 104, y from -84 to -4

\arrow <6pt> [.2,.67] from 2 -2 to 80 -80

\plot 32 -32 62 -2 /

 \plot 16 -16 30 -2 /

 \plot 48 -16 34 -2 /

 \plot 8 -8 14 -2 /

 \plot 24 -8 18 -2 /

 \plot 40 -8 46 -2 /

 \plot 56 -8 50 -2 /

 \plot 4 -4 6 -2 /

 \plot 12 -4 10 -2 /

 \plot 20 -4 22 -2 /

 \plot 28 -4 26 -2 /

 \plot 36 -4 38 -2 /

 \plot 44 -4 42 -2 /

 \plot 52 -4 54 -2 /

 \plot 60 -4 58 -2 /



 \plot 99 -29 64 -64 /

 \plot 66 -2 96 -32 /

 \plot 70 -2 68 -4 /

 \plot 74 -2 76 -4 /

 \plot 78 -2 72 -8 /

 \plot 82 -2 88 -8 /

 \plot 86 -2 84 -4 /

 \plot 90 -2 92 -4 /

 \plot 94 -2 80 -16 /


\setdots <3pt>
\putrule from -4.8 -4 to 102 -4
\putrule from -4.5 -8 to 102 -8
\putrule from -2 -16 to 102 -16
\putrule from -1.7 -32 to 102 -32
\putrule from -1.7 -64 to 102 -64

\put {$\vdots$} at 32 3
\put {$\vdots$} at 64 3

\put {$\dots$} [l] at 103 -6
\put {$\dots$} [l] at 103 -48

\put {$H_{-3}$} [l] at -12 -64
\put {$H_{-2}$} [l] at -12 -32
\put {$H_{-1}$} [l] at -12 -16
\put {$H_0$} [l] at -12 -8
\put {$H_1$} [l] at -12 -4
\put {$\vdots$} at -10 3
\put {$\vdots$} [B] at -10 -70
\put {$\omega$} at 82 -82

\put {\scriptsize $\bullet$} at 44 -4
\put {\scriptsize $1$} at 44 -6.5
\put {\scriptsize $0$} at 41 -12
\put {\scriptsize $1$} at 43 -24
\put {\scriptsize $0$} at 44 -48
\put {\scriptsize $0$} at 71.5 -75

\endpicture
$$

\vspace{.1cm}

\begin{center}
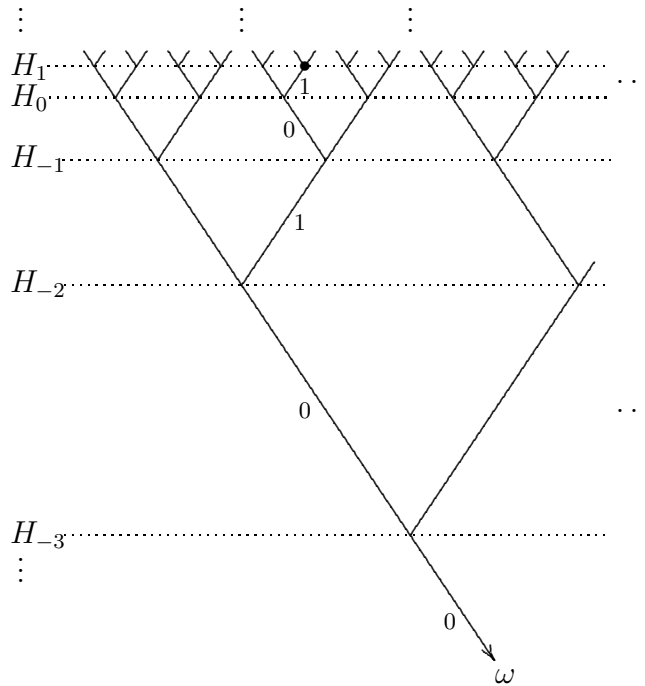
\centerline{\it Figure 1}
\end{center}


\begin{imp}\label{tree-seq} {\bf Tree and sequences.} \rm
Now consider the set $\Sigma_q$ of all sequences 
$\bigl( \sigma(n) \bigr)_{n \le 0}$ over $\Z_q$ with finite support
$\{ n : \sigma(n) \ne 0 \}$. We denote by $\tau$ the (negative) shift,
$\tau\sigma(n) = \sigma(n-1)$. Then the set $\Sigma_q \times \Z$
carries the structure of $\T_q$ in horocylic layers as above:
the $k$-th horocycle is $H_k = \Sigma_q \times \{k\}$, and
the predecessor of vertex $x = (\sigma, k)$ is $x^- = (\tau\sigma,k-1)$.
This corresponds to labelling the edges of $\T_q$ by elemets of $\Z_q$
such that all edges on the ray from $\om$ to $o$ have label $0$, and
for every vertex $x$ and every $\ell \in \Z_q$ there is a successor $y$
among the $q$ successors of $x$ such that the edge $[x,y]$ carries
label $\ell$. See Figure 1: the origin is the leftmost point on 
the horocycle $H_0$, 
we have indicated the labels on the edges that lead to the point marked with
a ``{\small $\bullet$}'', and that point has coordinates $(\sigma,k)$ with
$\sigma = (\dots,0,\dots,0,1,0,1)$ and $k=-1$. 
\end{imp}

Now consider two trees $\T_q$ and $\T_r$ with roots $o_1$ and $o_2$ and
reference ends $\om_1$ and $\om_2$, respectively.

\begin{dfn}\label{DLdef}
The Diestel-Leader graph $\DL(q,r)$ is
$$
\DL(q,r) = \{ x_1x_2 \in \T_q \times \T_r : \hor(x_1)+\hor(x_2) = 0 \}\,,
$$
and neighbourhood is given by
$$
x_1x_2 \sim y_1y_2 \iff x_1 \sim y_1 \AND x_2 \sim y_2\,.
$$
\end{dfn}

To visualize $\DL(q,r)$, draw $\T_q$ in horocyclic layers with
$\om_1$ at the top and $\bd^*\T_q$ at the bottom, and right to it $\T_r$
in the same way, but upside down, with the respective horocycles $H_k(\T_q)$
and $H_{-k}(\T_r)$ on the same level. Connect the two origins $o_1$, $o_2$ by
an elastic spring. It is allowed to move along each of the two trees, 
may expand infinitely, but must always remain in horizontal position. 
The vertex set of $\DL_{q,r}$ consists of all admissible positions of 
the spring. From a position
$x_1x_2$ with $\hor(x_1) + \hor(x_2) =0$ the spring may
move downwards to one of the $r$ successors of $x_2$ in $\T_r$, and at the same 
time to the predecessor of $x_1$ in $\T_q$, or it may move upwards 
in the analogous way. Such a move corresponds
to going to a neighbour of $x_1x_2$. Figure 2 depicts
$\DL(2,2)$.


$$
\beginpicture
\setcoordinatesystem units <3mm,3.5mm> 

\setplotarea x from -4 to 30, y from -3.8 to 6.4
\arrow <5pt> [.2,.67] from 4 4 to 1 7
\put{$\omega_1$} [rb] at 1.2 7.2

\put{$o_1$} [lb] at  8.15 0.2

\plot -4 -4       4 4         /         
\plot 4 4         12 -4          /      
\plot -2 -2       -2.95 -4 /            
\plot -.5 -2      -1.9 -4     /         
\plot -.5 -2      -.85 -4   /           
\plot 1 -2        .2 -4      /          
\plot 1 -2        1.25  -4     /        
\plot 2.5 -2      2.3  -4     /         
\plot 2.5 -2      3.35 -4      /        
\plot 5.5 -2      4.65  -4    /         
\plot 5.5  -2     5.7  -4      /        
\plot 7  -2       6.75   -4   /         
\plot 7 -2        7.8  -4      /        
\plot 8.5  -2     8.85  -4    /         
\plot 8.5  -2     9.9  -4      /        
\plot 10  -2      10.95  -4   /         
\plot 0  0        -.5  -2      /        
\plot 2 0         1 -2     /            
\plot 2 0         2.5   -2     /        
\plot 6 0         5.5    -2    /        
\plot 6 0         7 -2         /        
\plot 8 0         8.5 -2       /        
\plot 2 2         2 0         /         
\plot 6 2         6 0         /         

\arrow <5pt> [.2,.67] from 22 -4 to 25 -7
\put{$\omega_2$} [lt] at 25.2 -7.2

\put{$o_2$} [rt] at  17.95 -.2

\plot 14  4       22 -4       /         
\plot 22 -4         30 4         /      
\plot 16 2       15.05 4 /              
\plot 17.5 2     16.1  4     /          
\plot 17.5 2     17.15  4   /           
\plot 19  2       18.2  4      /        
\plot 19 2         19.25   4     /      
\plot 20.5 2       20.3   4     /       
\plot 20.5  2      21.35 4      /       
\plot 23.5 2       22.65  4    /        
\plot 23.5  2      23.7   4      /      
\plot 25   2       24.75   4   /        
\plot 25  2        25.8  4      /       
\plot 26.5   2     26.85  4    /        
\plot 26.5   2     27.9   4      /      
\plot 28   2      28.95   4   /         
\plot 18  0        17.5  2      /       
\plot 20 0         19  2     /          
\plot 20 0         20.5   2     /       
\plot 24 0         23.5    2    /       
\plot 24 0         25 2         /       
\plot 26 0         26.5 2       /       
\plot 20 -2        20 0         /       
\plot 24 -2        24 0         /       
\put {$\circ$} at 8 0
\put {$\circ$} at 18 0
\plot 8.25 0  12.1 0 /
\plot 13.9 0 17.78 0 /
\plot    12.1   0    12.25 .4    12.25 -.4   12.55 .4   12.55 -.4
        12.85 .4   12.85 -.4  13.15 .4   13.15 -.4  13.45 .4
        13.45 -.4  13.75 .4   13.75 -.4  13.9 0     13.9  0 /

\setdashes <2pt>
\putrule from -4.5 -7  to  12.5 -7
\putrule from  13.5 7  to  30.5 7

\put {$\bd^*\T_q$} [r] at -5 -7
\put {$\bd^*\T_r$} [l] at 31 7

\put {$\vdots$} at 4 -5.2
\put {$\vdots$} at 22 5.5

\endpicture
$$
\vspace{.4cm}

\begin{center}
{\it Figure 2}
\end{center}

\vspace{.4cm}

As the reference point in $\DL(q,r)$, we choose $o=o_1o_2$. We shall keep in
mind that $\T_q$ is the first and $\T_r$ the second tree; when $r=q$, it will be
sometimes convenient to write $\T^1$ and $\T^2$ for the first and second trees,
both copies of $\T_q$.

\vspace{.4cm}

Next, we explain what the lamplighter group $\Z_q \wr \Z$ has to do with 
$\DL(q,q)$. Let $(\eta,k) \in \Z_q \wr \Z$, and recall that $\eta: \Z \to \Z_q$ 
is a finitely supported configuration. We identify $(\eta,k)$ with the vertex 
$x_1x_2 \in \DL(q,q)$, where according to (\ref{tree-seq}), the vertices $x_i$ 
are given by 
\begin{equation}\label{identif}
\begin{gathered}
x_1 = (\eta_k^-,k) \AND x_2=(\eta_k^+,-k)\,,\quad\mbox{where}\\
\eta_k^- = \bigl(\eta(k+n)\bigr)_{n \le 0} \AND 
\eta_k^+ = \bigl(\eta(k+1-n)\bigr)_{n \le 0}\,,
\end{gathered}
\end{equation}
that is, we split $\eta$ at $k$, with $\eta_k^- = \eta|_{(-\infty\,,\,k]}$ 
and $\eta_k^+=\eta|_{[k+1\,,\,\infty)}$,
both written as sequences over the non-positive integers.

This is clearly a one-to-one correspondence between $\DL(q,q)$ and 
$\Z_q \wr \Z$, and it is also straighforward that this group acts 
transitively and fixed-point-freely on the graph: the action of $m \in \Z$ 
is given by $x_1x_2 = (\sigma_1,k)(\sigma_2,-k) 
\mapsto y_1y_2 = (\sigma_1,k+m)(\sigma_2,-k+m)$,
and the action of the group of configurations is pointwise addition modulo $q$ 
in the obvious way; the reader is invited to work out the simple details. We 
have to determine  the symmetric set of generators of our group with respect 
to which $\DL(q,q)$ is its Cayley graph. (Here, we mean the \emph{right} Cayley
graph, where an edge corresponds to multiplying with a generator on the right.)

Stepping from a vertex $x_1x_2$ to $y_1x_2^-$, where $y_1$ is one of 
the successors of $x_1 \in H_k(\T^1)$ (horocycle in the first tree)
means that the lamplighter walks from position $k$ to $k+1$ and then 
switches the lamp at the new position to some state in $\Z_q$. 
Thus, the ``downward'' edges of this 
type correspond to multiplying on the right with the group elements 
$(\de_1^{\ell},1)$, $\ell \in \Z_q$, where $\de_k^{\ell}$ is the 
configuration with value $\ell$ at $k$ and $0$ elsewhere.
On the other hand, we have the ``upward'' edges from $x_1x_2$ to $x_1^-y_2$, 
where  $y_2$ is one of the successors of $x_2 \in H_{-k}(\T^2)$ 
(horocycle in the second tree). They correspond to multiplying on the right with the inverses 
of the above generators,  i.e., the elements $(\de_0^{\ell},-1)$, where 
$\ell \in \Z_q$. 

Thus the simple random walk on $\DL(q,q)$ is the following lamplighter walk: 
its \emph{law,} the probability measure on $\Z_q \wr \Z$ that describes
the one step transition probabilites, is equidistribution on
\begin{equation}\label{walk-switch}
\{ (\de_1^{\ell},1)\,,\;(\de_0^{\ell},-1) : \ell \in \Z_q \}\,.
\end{equation}
If at some step, the lamplighter stands at $k \in \Z$, (s)he chooses with 
equal probability either to  step to $k+1$ and then to switch the lamp at 
$k+1$ to a random state, or (s)he chooses to
switch the lamp at $k$ to a random state (before leaving $k$) and then to 
step to $k-1$. While this is a symmetric random walk on $\Z_q \wr \Z$, 
resp. $\DL(q,q)$, this type of action does not appear ``symmetric'' from 
the point of view of the lamplighter. For this reason, other types of 
``simple'' random walks have been considered in the past: the one whose 
law is equidistribution on
\begin{equation}\label{walk-or-switch}
\{ ({\mathbf 0},\pm 1)\,,\;(\de_0^{\ell},0) : 0 \ne \ell \in \Z_q \}
\end{equation}
(``walk or switch''), and the one where the lamplighter standing at $k$ 
first switches the lamp where he stands to a random state, then walks to
$k \pm 1$, and then switches the lamp at the arrival point to a random state
(``switch-walk-switch''). 
The corresponding law is equidistribution on
\begin{equation}\label{switch-walk-switch}
\{ (\de_0^{\ell}+\de_{\pm 1}^{m},\pm 1) : \ell,m \in \Z_q \}
\end{equation}

Harmonic functions for the ``walk or switch'' model cannot be 
determined by the methods that we elaborate here, since it is not 
very well adapted to the structure of $\DL(q,q)$, see the comments at the end.

On the other hand, the ``switch-walk-switch'' model (\ref{switch-walk-switch})
corresponds to simple random walk on the following modification of 
$\DL(q,q)\,$: in the first of the two trees, we add edges between
every vertex and the siblings of its predecessor (i.e., its ``uncles''),
and the resulting neighbourhood relation on the horocyclic product is
as in Definition \ref{DLdef}. It will be easy to adapt our results to
this random walk. 

In the first place, we shall study the following slight generalization 
$P = P_{\al}$ of simple random walk on $\DL(q,r)$, where $0 < \al < 1$. 
For $x_1x_2 \in \DL(q,r)$
\begin{equation}\label{random-walk}
p(x_1x_2,y_1y_2) = \begin{cases} \al/q 
                   & \text{if}\; y_1^- = x_1 \;\text{and}\;y_2=x_2^-\\
                                 (1-\al)/r 
		   & \text{if}\; y_1 = x_1^- \;\text{and}\;y_2^-=x_2\\
                 0 & \text{otherwise.}
                   \end{cases}
\end{equation}
$P$ acts on functions $h: \DL(q,r) \to \R$ by
$$
Ph(x_1x_2) = \sum_{y_1y_2}  p(x_1x_2,y_1y_2)h(x_1x_2)
$$  
A \emph{harmonic,} or more precisely, \emph{$P_{\alpha}$-harmonic  function,}  is 
one that satisfies $Ph = h$.

We can consider the projections $P_1 = P_{1,\al}$ and $P_2 = P_{2,1-\al}$
of $P_{\al}$ on  $\T_q$ and $\T_r$, respectively:
\begin{equation}\label{projections}
p_1(x_1,y_1) = \begin{cases} \al/q & \text{if}\; y_1^- = x_1 \\
                             (1-\al) & \text{if}\; y_1 = x_1^- \\
			     0 & \text{otherwise,}
                                    \end{cases}
\qquad
p_2(x_1,y_2) = \begin{cases} \al & \text{if}\; y_2 = x_2^- \\
                             (1-\al)/r & \text{if}\; y_2^- = x_2 \\
			     0 & \text{otherwise.}
                                    \end{cases}
\end{equation}

The following is straightforward.

\begin{lem}\label{lift-harmonic}
{\rm (a)} If $h_1$ is a $P_1$-harmonic function on $\T_q$, then 
$h(x_1x_2) = h_1(x_1)$, $x_1x_2 \in \DL(q,r)$, defines a $P$-harmonic function
on $\DL(q,r)$.\\

{\rm (b)} If $h_2$ is a $P_2$-harmonic function on $\T_q$, then 
$h(x_1x_2) = h_2(x_2)$, $x_1x_2 \in \DL(q,r)$, defines a $P$-harmonic function
on $\DL(q,r)$.\\

\end{lem}

Our first main result is the following.

\begin{thm}\label{split-theorem}
If $h$ is a non-negative $P$-harmonic function on $\DL(q,r)$, then there
are non-negative $P_i$-harmonic functions $h_i\,$, $i=1,2$, on\/ $\T_q$ and\/
$\T_r$, respectively, such that
$$
h(x_1x_2) = h_1(x_1) +  h_2(x_2) \quad \text{for all}\; x_1x_2 \in \DL(q,r)
$$
\end{thm}

Conversely, it is of course clear that every sum of the latter form defines
a $P$-harmonic function. Recall that in this type of decomposition,  $x_2$ 
cannot vary independently of $x_1$, since one must have 
$\hor(x_1) + \hor(x_2) = 0$. 

The next short section contains some basic preparatory material 
for the proof of Theorem \ref{split-theorem}.

\medskip

\section{Basic results about harmonic functions}\label{basics}

Let $X$ be a denumerable set and $P = \bigl(p(x,y)\bigr)_{x,y \in X}$
the stochastic transition matrix of a Markov chain $(Z_n)_{n \ge 0}$ on $X$. 
%
%
We write $\Pr_x$ for probability conditioned to the starting point
$Z_0=x$. The $n$-step transition probability  $p^{(n)}(x,y) = \Pr_x[Z_n=y]$
is the $(x,y)$-entry of the matrix power $P^n$. We assume that 
$P$ is \emph{irreducible:} $\forall\ x,y\ \exists\ n: p^{(n)}(x,y) > 0$.

As above, a function $h$ on $X$ is called \emph{$P$-harmonic} or just
\emph{harmonic at} $x$, if
$Ph(x) = h(x)$, where $Ph(x) = \sum_y p(x,y)h(y)$. It is called harmonic
when it is harmonic at each $x$. 

For a subset $A \subset X$, we define the stopping time
$$
s^A = \inf \{ n \ge 0 : Z_n \in A \}\,.
$$
For $y \in X$, we write $s^y = s^{\{y\}}$. Given $x,y \in X$, let
\begin{equation}\label{hit-def}
F(x,y) = {\Pr}_x[s^y < \infty] \AND F^A(x,y) = {\Pr}_x[s^y \le s^A\,,\;s^y < \infty]
\end{equation}
Thus, $F(x,y)$ is the probability to ever reach $y$, starting from $x$.
The function $F(\cdot,y)$ is harmonic in $X \setminus \{y\}$.
Furthermore, if $y \in A$, then
$F^A(\cdot,y)$ is harmonic in $X \setminus A$. 

\begin{imp}\label{finite} {\bf Harmonic functions on finite sets.} 
\rm Let $S$ be a finite subset of $X$. Define its \emph{boundary}
and \emph{interior} by
$$
\bd S = \{ y \in S : p(y, X \setminus S) > 0\} \AND S^o = S \setminus \bd S\,.
$$
For the sake of simplicity, we assume that the restriction of $P$ to $S^o$ is
irreducible. We define
$$
\HH(P,S) = \{ h: S \to \R \mid h \;\text{is harmonic in}\; S^o \}
$$
The following is very well known.

\begin{pro}\label{dirichlet} Under the above assumptions,
the functions $F^{\bd S}(\cdot, y)\,$, $y \in \bd S$, constitute
a basis of the linear space $\HH(P,S)$. Every $h \in \HH(P,S)$ is uniquely
respresented as
$$
h = \sum_{y \in \bd S} F^{\bd S}(\cdot,y)\,h(y).
$$
\end{pro}

\begin{proof} The functions $F^{\bd S}(\cdot, y)\,$, $y \in \bd S$,
are linearly independent, since $F^{\bd S}(x, y)=\de_x(y)$ for $x, y \in \bd S$.
Given $h \in \HH(P,S)$, let $g = \sum_{y \in \bd S} F^{\bd S}(\cdot,y)\,h(y)$.
Then $g \in \HH(P,S)$, and $(g-h)|_{\bd S} \equiv 0$. By the 
\emph{Minimum Principle,} every function in $\HH(P,S)$ attains its
minimum (and its maximum) on the boundary. Therefore $g = h$ on $S$. 
\end{proof}

\end{imp}

\begin{imp}\label{positive-harmonic} {\bf Positive and 
minimal harmonic functions.} \rm Regarding the following material,
see {\sc Woess} \cite{Wbook}, \S 24 for a more detailed outline and
many references. 

We return to the infinite set $X$ with irreducible transition matrix $P$. 
For the sake of simplicity, we assume that $P$ has \emph{finite range,}
i.e., $\{ y: p(x,y) > 0 \}$ is finite for all $x \in X$. 
The set $\HH^+ = \HH^+(P,X)$ of 
non-negative $P$-harmonic functions constitutes a convex cone that is closed
in the topology of pointwise convergence. We choose a reference point
$o \in X$. Then the set $\BB = \{ h \in \HH^+ : h(o)=1 \}$
is a compact, convex base of the cone $\HH^+$. Its extremal elements are called
\emph{minimal harmonic functions.} Thus, $h \in \HH^+$ is minimal 
if
$$
h(o) = 1 \AND h \ge h_1 \in \HH^+ \Longrightarrow h_1/h \equiv \text{constant.}
$$
The set $\BB_{\min}$ of minimal harmonic functions is a Borel subset of $\BB$,
and every $h \in \HH^+$ is an integral of minimal ones with respect to
a Borel measure on $\BB_{\min}$. This can be made more precise by the following
construction. Define the \emph{Martin kernel}
$$
K(x,y) = F(x,y)/F(o,y)\,.
$$
The \emph{Martin compactification} is the smallest metrizable compactification
of $X$  containing $X$ as a discrete, dense subset, and to which
all functions $K(x,\cdot)$, $x \in X$, extend continuously. The 
\emph{Martin boundary} $\MM = \MM(P)$ is the ideal boundary added to $X$
in this compactification. Then every minimal harmonic function is
of the form $K(\cdot,\xi)$ for some $\xi \in \MM$, and the set
$$
\MM_{\min} = \{ \xi \in \MM : K(\cdot,\xi) \;\text{is minimal harmonic} \}
$$
is a Borel set. The \emph{Poisson-Martin Representation Theorem} says
that for every $h \in \HH^+$ there is a unique Borel measure $\nu^h$ on $\MM$
with $\nu^h(\MM \setminus \MM_{\min})=0$ such that
$$
h(x) = \int_{\MM} K(x,\cdot)\,d\nu^h \quad \forall\ x \in X\,.
$$
Furthermore, considering the constant harmonic function $\mathbf 1$, 
we set $\nu = \nu^{\mathbf 1}$. Then every \emph{bounded} harmonic 
function $h$ has a unique representation as above, where 
$d\nu^h(\xi) = \varphi(\xi)\,d\nu(\xi)$ with $\varphi \in L^{\infty}(\MM,\nu)$.
The probability space $(\MM,\nu)$ is a model of the \emph{Poisson boundary}
of the random walk. While the Martin boundary is a topological object,
the Poisson boundary is a measure theoretical one, and finding it means
to determine it up to isomorphisms between measure spaces. 

See {\sc Kaimanovich and Vershik} \cite{KaiVer} for a profound introduction and 
impressive results regarding Poisson boundaries of random walks on groups,
{\sc Kaimanovich} \cite{Kai} for lamplighter groups over $\Z^d$ and other
semidirect products, and {\sc Kaimanovich and Woess} \cite{KaiWoe} for 
Poisson boundaries of random walks on homogenenous graphs, including the 
Diestel-Leader graphs.
\end{imp}

\begin{imp}\label{trees} {\bf Harmonic functions on trees.} \rm
Next, let us suppose that $X=T$ carries the structure of an infinite,
locally finite tree. We assume that $P$ is of \emph{nearest neighbour}
type, i.e., 
\begin{equation}
p(x,y) > 0 \iff x \sim y \; \text{in}\; T\,.
\end{equation}
($\sim$ denotes neighbourhood.) We also assume that the random walk
(Markov chain) with transition matrix $P$ is \emph{transient,} that is,
$\sum_n p^{(n)}(x,y) < \infty$ for some ($\!\!\iff\!$ all) $x,y \in T$.
Geodesics and boundary of $T$ are defined as in \S \ref{geometry},
with the general tree $T$ in the place of $\T_q$. The results regarding
the Martin compactification in this setting are contained in the
seminal paper by {\sc Cartier} \cite{Car}.

The basic link between tree structure and random walk is the following 
well-known lemma, see e.g. {\sc Cartier} \cite{Car}, or {\sc Woess} \cite{Wbook}, 
Lemmas 1.23 and 1.13(d).

\begin{lem}\label{tree-lemma}
For a nearest neighbour random walk on a tree $T$, 
$$
F(x,y) = F(x,w)F(w,z) 
\quad \text{for all}\;x,y \in T\;\text{and}\; w \in \geo{x\,y}\,.
$$
Furthermore, if $x \sim y$, then 
$$
F(y,x) = p(y,x) + \sum_{w \ne x} p(y,w)F(w,y)F(y,x)\,.
$$
\end{lem}

For $x, y \in T$, let $c= x \wedge y$ be their confluent with respect to $o$.
Then the lemma implies that $K(x,y) = K(x,c)$. From here, the following
is almost immediate.

\begin{pro}\label{tree-martin}
Suppose that $P$ defines a transient nearest neighbour random walk on
the locally finite tree $T$ with root $o$. Then the Martin compactification
is the end compactification $\wh T$, and for $\xi \in \bd T$, the 
Martin kernel is given by
$$
K(x,\xi) = K(x,c) = \frac{F(x,c)}{F(o,c)}\,,\quad 
\text{where}\quad c = x \wedge \xi\,.
$$
Furthermore, each function $K(\cdot,\xi)$, $\xi \in \bd T$, 
is minimal harmonic.
\end{pro}

For various different proofs, see {\sc Cartier} \cite{Car}, 
{\sc Picardello, Taibleson and Woess} \cite{PicTaiWoe}, or {\sc Woess} 
\cite{Wbook}, \S 26, or also the one which is implicit in the proof of 
Theorem \ref{split-theorem} below.   
\end{imp}

\begin{exa}\label{alpha-walk} Consider the random walk on $\T_q$ with transition matrix
$P_1=P_{1,\al}$, defined in \ref{projections}.
It is clear that for this random walk, the probabilities
$$
F_1^- = F_1(x,x^-) \AND  F_1^+ = F_1(x^-,x)
$$
are independent of $x \in \T_q$ ($x^-$ is the predecessor with respect to 
$\om$). Using Lemma \ref{tree-lemma}, we find the two quadratic equations
$$
F_1^- = (1-\al) + \al (F_1^-)^2 \AND 
F_1^+ = \frac{\al}{q} + (q-1)\frac{\al}{q}F_1^-F_1^+ + (1-\al)(F_1^+)^2\,.
$$
Among the two solutions of each equation, the smaller one is the right one (compare e.g.
with the generating functions argument in the proof of Lemma 1.24 in 
{\sc Woess} \cite{Wbook}).
Thus
\begin{equation}\label{F-computation}
F_1^- = \begin{cases} \dps \frac{1-\al}{\al} 
                           &\dps \text{if}\; \al \ge \frac12\,,   \\[7pt]
                           \dps 1& \dps \text{if}\; \al \le \frac12\,, 
                                    \end{cases} 
\qquad
F_1^+ = \begin{cases} \dps \frac{1}{q} 
                           & \dps \text{if}\; \al \ge \frac12\,,    \\[7pt]
                           \dps   \frac{\al}{(1-\al)q} 
			   & \dps \text{if}\; \al \le \frac12\,.
                                    \end{cases}
\end{equation}
We can now compute the associated Martin kernels $K_1(\cdot,\xi)$, $\xi \in \bd \T_q$.
First, since $x \wedge \om = x \cf o$ (where $\wedge$ and $\cf$ denote confluents with 
respect to $o$ and $\om$), it is immediate that
\begin{equation}\label{martin-omega}
K_1(x,\om) = (F_1^-)^{\hor(x)} = \begin{cases}
                                      \dps \left(\frac{1-\al}{\al}\right)^{\mbox{\footnotesize $\hor(x)$}}
                                      & \dps \text{if}\; \al \ge \frac12\,,  \\[7pt]
                                       \dps 1 & \dps \text{if}\; \al \le \frac12\,.
                                    \end{cases}
\end{equation}
Next, if $\xi \in \bd^*\T_q$, we set $c = x \wedge \xi$ and write
$k = d(o,o \cf x)$, $l = d(x,o \cf x)$, so that $\hor(x) = l-k$. 
We distinguish two cases, see Figure 3.
$$
\beginpicture 

\setcoordinatesystem units <1mm,1.6mm>

\setplotarea x from -2 to 100, y from -10 to 10


\arrow <6pt> [.2,.67] from 0 0 to 17 -17
\arrow <6pt> [.2,.67] from 6 -6 to 16.5 4.5
\plot 12 -12   28 4 /

\put {$\omega$} at 18 -18
\put {$\xi$} at 18 6
\put {$o$} at -1.5 1.5
\put {$x$} at 29.5 5.5
\put {$c$} at 4 -7.5
\put {$o \cf x$}[r] at 10.5 -13
\multiput {\scriptsize $\bullet$} at 0 0  6 -6  12 -12  28 4 /


\arrow <6pt> [.2,.67] from 70 0 to 87 -17
\arrow <6pt> [.2,.67] from 88 -6 to 77.5 5.5
\plot 82 -12   98 4 /

\put {$\omega$} at 88 -18
\put {$\xi$} at 76 6
\put {$o$} at 68.5 1.5
\put {$x$} at 99.5 5.5
\put {$c$} at 89.5 -7.5
\put {$o \cf x$}[r] at 80.5 -13
\multiput {\scriptsize $\bullet$} at 70 0  82 -12  88 -6  98 4 /

\put{\it Figure 3}[c] at 49 -22
\endpicture
$$

\vspace{.4cm}



\noindent \emph{Case 1.} $c$ lies between $o$ and $o \cf x$. Let $s = d(o,c)$. Then
$$
K_1(x,\xi) = \frac{(F_1^-)^l(F_1^+)^{k-s}}{(F_1^-)^s} = K_1(x,\om)(F_1^-F_1^+)^{k-s}
= K_1(x,\om)(F_1^-F_1^+)^{\mbox{\footnotesize $\hor(o\cf\xi)-\hor(x\cf\xi)$}}
$$
\emph{Case 2.} $c$ lies between $o \cf x$ and $x$. Let $r = d(x,c)$. Then
$$
K_1(x,\xi) = \frac{(F_1^-)^r}{(F_1^-)^k(F_1^+)^{l-r}} = K_1(x,\om)(F_1^-F_1^+)^{r-l}
= K_1(x,\om)(F_1^-F_1^+)^{\mbox{\footnotesize $\hor(o\cf\xi)-\hor(x\cf\xi)$}}
$$
We write $\hor(x,\xi) = d(x,c) - d(o,c)$, the horocycle number with respect to $\xi$, while
$\hor(x) = \hor(x,\om)$. Also, we set $\rho =  (F_1^-F_1^+)^{1/2}$. Then we find in both cases 
\begin{equation}\label{martin-xi}
K_1(x,\xi) = K_1(x,\om)\, \rho^{\mbox{\footnotesize $\hor(x,\xi) - \hor(x)$}}
\,,\quad\text{where}\quad
\rho = \min \left\{ \frac{1-\al}{\al q}\,, \frac{\al}{(1-\al)q} \right\}^{1/2}\,.
\end{equation}
In particular, if $\al =1/2$ then 
\begin{equation}\label{martin-1}
K_1(x,\om) = 1 \AND 
K_1(x,\xi) = q^{\mbox{\footnotesize $(\hor(x) - \hor(x,\xi))/2$}}\quad 
\text{for}\quad \xi \in \bd^*\T_q\,.
\end{equation}
\end{exa}

\medskip

\section{Minimal harmonic functions on $\DL(q,r)$}\label{principal}

After all these preliminaries, the proof of Theorem \ref{split-theorem} depends in the first place 
on the way how we look at the underlying structure.

\begin{proof}[Proof of Theorem \ref{split-theorem}]
In $\DL(q,r)$, consider the subgraph spanned by all vertices $x_1x_2$ with 
$-n \le \hor(x_1) \le n$. It is not connected. We denote by $S = S^{(n)}$ the 
connected component of the root $o_1o_2$. 

Let $a_1 = a_1^{(n)} \in \T_q$ be the vertex on $\geo{o_1\,\om_1}$ at distance
$d(a_1,o_1) = n$. Then $a_1$ can be viewed as the root of the $q$-ary 
rooted tree 
$S_1 = S_1^{(n)} = \{ x_1 \in \T_q : -n \le \hor(x_1) \le n\,,\; a_1 \in \geo{x_1\,\om_1}\,\}$
of height $2n$, whose set of leaves (elements with $\hor(x_1)=n$) is denoted
$\bd^* S_1$. Analogously, we define $a_2 = a_2^{(n)}$, the $r$-ary rooted tree 
$S_2 = S_2^{(n)}$, and it set of leaves $\bd^*S_2$. Then 
$$
S = \{ x_1x_2 \in S_1 \times S_2 : \hor(x_1)+\hor(x_2)=0\}
$$
is the horocyclic product of $S_1$ and $S_2$, see Figure 4. 

$$
\beginpicture
\setcoordinatesystem units <3mm,3.5mm> 

\setplotarea x from -4 to 30, y from -2.3 to 2.3
\put{$o_1$} [lb] at  8.15 0.2

\plot -4 -4       4 4         /         
\plot 4 4         12 -4          /      
\plot -2 -2       -2.95 -4 /            
\plot -.5 -2      -1.9 -4     /         
\plot -.5 -2      -.85 -4   /           
\plot 1 -2        .2 -4      /          
\plot 1 -2        1.25  -4     /        
\plot 2.5 -2      2.3  -4     /         
\plot 2.5 -2      3.35 -4      /        
\plot 5.5 -2      4.65  -4    /         
\plot 5.5  -2     5.7  -4      /        
\plot 7  -2       6.75   -4   /         
\plot 7 -2        7.8  -4      /        
\plot 8.5  -2     8.85  -4    /         
\plot 8.5  -2     9.9  -4      /        
\plot 10  -2      10.95  -4   /         
\plot 0  0        -.5  -2      /        
\plot 2 0         1 -2     /            
\plot 2 0         2.5   -2     /        
\plot 6 0         5.5    -2    /        
\plot 6 0         7 -2         /        
\plot 8 0         8.5 -2       /        
\plot 2 2         2 0         /         
\plot 6 2         6 0         /         

\multiput {\scriptsize $\bullet$} at   4 4   
       -4 -4    -2.95 -4    -1.9 -4    -.85 -4    .2 -4     1.25   -4   2.3  -4     3.35 -4    
        4.65  -4    5.7  -4     6.75  -4    7.8  -4      8.85  -4    9.9  -4      10.95  -4   12 -4     /

\put{$a_1$}[b] at 4 4.5
\put{$\bd^*S_1$}[t] at 4 -5


\put{$o_2$} [rt] at  17.95 -.2

\plot 14  4       22 -4       /         
\plot 22 -4         30 4         /      
\plot 16 2       15.05 4 /              
\plot 17.5 2     16.1  4     /          
\plot 17.5 2     17.15  4   /           
\plot 19  2       18.2  4      /        
\plot 19 2         19.25   4     /      
\plot 20.5 2       20.3   4     /       
\plot 20.5  2      21.35 4      /       
\plot 23.5 2       22.65  4    /        
\plot 23.5  2      23.7   4      /      
\plot 25   2       24.75   4   /        
\plot 25  2        25.8  4      /       
\plot 26.5   2     26.85  4    /        
\plot 26.5   2     27.9   4      /      
\plot 28   2      28.95   4   /         
\plot 18  0        17.5  2      /       
\plot 20 0         19  2     /          
\plot 20 0         20.5   2     /       
\plot 24 0         23.5    2    /       
\plot 24 0         25 2         /       
\plot 26 0         26.5 2       /       
\plot 20 -2        20 0         /       
\plot 24 -2        24 0         /       

\multiput {\scriptsize $\bullet$} at   22 -4       
      14  4     15.05 4   16.1  4     17.15  4    18.2  4     19.25   4   20.3   4     21.35 4   
      22.65  4    23.7   4     24.75  4    25.8  4    26.85  4    27.9   4    28.95   4   30 4         /     

\put{$a_2$}[t] at 22  -4.5
\put{$\bd^*S_2$}[b] at 22  5

\put {$\circ$} at 8 0
\put {$\circ$} at 18 0
\plot 8.25 0  12.1 0 /
\plot 13.9 0 17.78 0 /
\plot    12.1   0    12.25 .4    12.25 -.4   12.55 .4   12.55 -.4
        12.85 .4   12.85 -.4  13.15 .4   13.15 -.4  13.45 .4
        13.45 -.4  13.75 .4   13.75 -.4  13.9 0     13.9  0 /

\endpicture
$$

\begin{center}
{\it Figure 4}
\end{center}

\vspace{.4cm}
One may imagine $S$ as a tetrahedron. Two of its faces are copies of $S_1$ that meet
at the common bottom edge $\bd^*S_1 \times \{a_2\}$, and the other two faces are 
copies of $S_2$ that meet at the common top edge $\{a_1\}\times \bd^*S_2$.
The boundary of $S_i$ is $\{a_i\} \cup \bd^* S_i$.

We now restrict $P$ to $S$, and also the projections $P_1$ to $S_1$ and $P_2$ to $S_2$.
Then the boundary of $S$ in the sense of (\ref{finite}) is
$$
\bd S = \bigl(\bd^*S_1 \times \{a_2\}\bigr) \cup \bigl(\{a_1\}\times \bd^*S_2\bigr)\,.
$$
As in  Lemma \ref{lift-harmonic}, if $h_i \in \HH(P_i,S_i)$, then it lifts to a function
in $\HH(P,S)$.

In particular, if $y_1 \in \bd^*S_1$ then $h(x_1x_2)= F_1^{\bd S_1}(x_1,y_1)$ defines  a function
in $\HH(P,S)$ with value one at $y_1a_2$ and value $0$ in $\bd S \setminus \{ y_1a_2 \}$.
But, by Proposition \ref{dirichlet}, these properties characterize the function
$x_1x_2 \mapsto F^{\bd S}(x_1x_2,y_1a_2)$ on $S$. Therefore, for all $x_1x_2 \in S$,
\begin{equation}\label{hitting}
\begin{aligned}
F^{\bd S}(x_1x_2,y_1a_2) &= F_1^{\bd S_1}(x_1,y_1) \quad 
\forall\ y_1 \in \bd^*S_1      
                          \quad\text{and}\\
F^{\bd S}(x_1x_2,a_1y_2) &= F_2^{\bd S_2}(x_2,y_2) \quad 
\forall\ y_2 \in \bd^*S_2 
\end{aligned}
\end{equation}
Applying Proposition \ref{dirichlet} once more, we see that every 
$h \in \HH(P,S)$ can be written uniquely as
\begin{equation}\label{finite-split}
\begin{aligned}
h(x_1x_2) &= h_1(x_1) + h_2(x_2) \quad \forall\ x_1x_2 \in S\,,\quad\text{where}\\[10pt]
h_1(x_1) &= \sum_{y_1 \in \bd^*S_1} F_1^{\bd S_1}(x_1,y_1) h(y_1a_2)
\AND\\
h_2(x_2) &= \sum_{y_2 \in \bd^*S_2} F_2^{\bd S_2}(x_2,y_2) h(a_1y_2)\,.
\end{aligned}
\end{equation}
This is true, in particular, if $h$ is $P$-harmonic on the whole of $\DL(q,r)$,
since its restriction to $S = S^{(n)}$ is in $\HH(P,S)$. Furthermore, if $h$ is non-negative
then so are $h_1$ and $h_2$. Note, however, that $h_i = h_i^{(n)}$ ($i=1,2$) \emph{depend on
$n$,} and it is by no means true that the restriction of $h_i^{(n+1)}$ to $S_i^{(n)}$
might coincide with $h_i^{(n)}$. We have to study the behaviour of $h_i^{(n)}$ when
$n \to \infty$, and this is the point where the assumption of non-negativity of $h$ will be used. 
We define 
$$
K_i^{(n)}(x_i,y_i) = \frac{F_i^{\bd S_i}(x_i,y_i)}{F_i^{\bd S_i}(o_i,y_i)}\,,\;
x_i \in S_i\,,\; y_i \in \bd S_i\,,\; S_i=S_i^{(n)}\,,\; i=1,2\,.
$$
Then we can rewrite the functions $h_i$ of (\ref{finite-split}) as
\begin{eqnarray*}
&&h_i(x_i) = \sum_{y_i \in \bd S_i} K_i^{(n)}(x_i,y_i) \la_i(y_i)\,,
\quad\text{where}\\
\la_1(a_1) \!\! &=&\!\! 0\,,\quad
\la_1(y_1) = h(y_1a_2)/F_1^{\bd S_1}(o_1,y_1)\,,\;y_1 \in \bd^*S_1\,,\\
\la_2(a_2) \!\! &=&\!\! 0\,, \quad
\la_2(y_2) = h(a_1y_2)/F_2^{\bd S_2}(o_2,y_2)\,,\;y_2 \in \bd^*S_2\,.
\end{eqnarray*}
Of course, also $\la_i(y_i) = \la_i^{(n)}(y_i)$ depends on $n$. 
For $\xi_1 \in \bd \T_q$, we define
$$
K_1^{(n)}(x_1,\xi_1) =  K_1^{(n)}(x_1,y_1)\,,\quad\text{if}\; \xi_1\in 
\T_q(o_1,y_1) \;\text{with}\;y_1 \in \bd S_1\,.
$$
Then   $K_1^{(n)}(x_1,\cdot)$ is locally constant (whence continuous) on 
$\bd \T_q$. (Recall that $\bd\T_q$ is compact and totally disconnected.)
We  define a non-negative Borel-measure $\nu_1^{(n)}$ on $\bd\T_q$
by
\begin{eqnarray*}
\nu_1^{(n)}\bigl(\bd\T_q(o_1,a_1^-)\bigr) &=& 0 \AND \\ 
\nu_1^{(n)}\bigl(\bd\T_q(o_1,w_1)\bigr)  &=& \la_1^{(n)}(y_1)\, q^{\hor(w_1)-n} 
\quad\text{if}\;
w_1 \in \T_q(o_1,y_1)\,,\;y_1 \in \bd^*S_1^{(n)}\,.
\end{eqnarray*}
This defines a finitely additive, non-negative measure on the 
semiring of all sets $\bd\T_q(o_1,x_1)$, $x_1\in \T_q \setminus \{o_1\}$. 
Since all these sets are open and compact, the measure is sigma-additive 
on that semiring and extends to a unique non-negative Borel measure on 
$\bd\T_q$. We proceed in precisely the same way on the second tree,
and also get a non-negative Borel measure $\nu_2^{(n)}$ on $\bd\T_r$, 
such that for all $x_1x_2 \in S^{(n)}\,$,
\begin{equation}\label{finite-integral}
h_1^{(n)}(x_1) = \int_{\bd\T_q} K_1^{(n)}(x_1,\cdot)\,d\nu_1^{(n)}  \AND
h_2^{(n)}(x_2) = \int_{\bd\T_r} K_2^{(n)}(x_2,\cdot)\,d\nu_2^{(n)} \,.
\end{equation}
Since $K_i^{(n)}(o_i,\cdot) \equiv 1$, we have 
$\nu_1^{(n)}(\bd\T_q) + \nu_2^{(n)}(\bd\T_r)  = h(o_1o_2)$ for all $n$.
Thus, by compactness (Helly's theorem), there are a subsequence $(n')$ 
 and non-negative measures $\nu_1$ on $\bd\T_q$ and $\nu_2$ on $\bd\T_r$
such that $\nu_i^{(n')} \to \nu_i$ weakly for $i=1,2$. 

If $x_1x_2 \in \DL(q,r)$, then we choose $n_0=n_0(x_1x_2)$ large enough 
such that for all $n \ge n_0$, the geodesics $\geo{o_i\,x_i}$ are contained 
in the interior of $S_i^{(n)}$, $i=1,2$. 

For every $\xi_i \in \bd\T_q$, resp. $\in \bd\T_r$, the confluent 
$c_i = x_i \wedge \xi_i$ is one of the finitely many points on 
$\geo{o_i\,x_i}$. It is clear that 
$F_i^{\bd S_i^{(n)}}(x_i,c_i) \to F_i(x_i,c_i)$ as $n \to \infty$, and the 
same is true with $o_i$ in the place of $x_i$. Therefore, using 
Lemma \ref{tree-lemma}, we find 
$$
K_i^{(n)}(x_i,\xi_i) = K_i^{(n)}(x_i,c_i) \to K_i(x_i,c_i) = K_i(x_i,\xi_i)\,,
$$ 
as $n \to \infty$, where $K_i(\cdot,\cdot)$ is the Martin kernel of 
$P_i$ on $\T_q$, resp. $\T_r$, $i=1,2$. Thus, 
$K_i^{(n)}(x_i,\cdot) \to K_i(x_i,\cdot)$ uniformly, and using 
(\ref{finite-integral}), we get
$$
h_1^{(n')}(x_1) \to \int_{\bd\T_q} K_1(x_1,\cdot)\,d\nu_1 =: h_1(x_1)\AND
h_2^{(n')}(x_2) \to \int_{\bd\T_r} K_1(x_2,\cdot)\,d\nu_2 =: h_2(x_2)\,.
$$
Then $h_1$ is $P_1$-harmonic on $\T_q$ and $h_2$ is
$P_2$-harmonic on $\T_r$, and $h(x_1x_2) = h_1(x_1)+h_2(x_2)$ for all 
$x_1x_2 \in \DL(q,r)$.
 \end{proof}

We remark that the last part is   is the key point of the argument, 
namely a way of recovering a Poisson integration formula
from the finite approximations $K_i^{(n)}(\cdot,\cdot)$ of the Martin kernel.

\begin{thm}\label{minimal-theorem}
{\rm (a)} Each of the functions 
$$
x_1x_2 \mapsto K_1(x_1,\xi_1)\,,\; \xi_1 \in \bd^*\T_q\,, \AND
x_1x_2 \mapsto K_2(x_2,\xi_2)\,,\; \xi_2 \in \bd^*\T_r\,,
$$
is minimal $P_{\al}$-harmonic on $\DL(p,q).$

\vspace{.2cm}

{\rm (b)} If $\al \ne 1/2$ then these are all minimal harmonic
functions.

\vspace{.2cm}

{\rm (c)} If $\al=1/2$, then these together with the constant function
$\mathbf 1$ are all minimal harmonic functions.
\end{thm}

\begin{proof} 
(a) Let $\xi_1 \in \bd^* \T_q$ and suppose that $K_1(x_1,\xi_1) \ge h(x_1x_2)$
for all $x_1x_2 \in \DL(q,r)$, where $h \ge 0$ is $P_{\al}$-harmonic.
By Theorem \ref{split-theorem}, $h(x_1x_2) = h_1(x_1) + h_2(x_2)$,
where $h_i \ge 0$ is $P_i$-harmonic, $i=1,2$. Then 
$K_1(\cdot,\xi_1) \ge h_1$. By Proposition \ref{tree-martin}, 
$h_1 = c\cdot K_1(\cdot,\xi_1)$, where $0 \le c \le 1$. 
If $c=1$ then we are done. Otherwise, 
$K_1(x_1,\xi_1) \ge c\cdot K_1(\cdot,\xi_1)+ h_2(x_2)$, that is,
$$
K_1(x_1,\xi_1) \ge \frac{1}{1-c}h_2(x_2) 
= \int_{\bd T_r} K_2(x_2,\cdot)\, d\nu_2
\quad \forall\ x_1x_2 \in \DL(q,r)\,,
$$
where $\nu_2$ is a non-negative Borel measure on $\bd \T^r$.

Setting $x_2 = o_2$, we obtain
$$
K_1(x_1,\xi_1) \ge \nu_2(\bd \T_r) \quad \forall\ x_1 \in H_0(\T_q)\,.
$$
If $x_1 \to \om_1$ then $\hor(x_1,\xi_1) \to \infty$.
Therefore, (\ref{martin-xi}) yields
$$
\nu_2(\bd \T_r) \le \lim_{x_1 \to \om_1,\, \hor(x_1)=0} 
K_1(x_1,\xi_1) = 0\,,
$$
and $\nu_2(\bd\T_r) = 0$. Therefore $h_2 \equiv 0$. 

This proves minimality of $K_1(\cdot,\xi_1)$ for all $\xi_1 \in \bd^*\T_q$.
Exchanging the roles of the two trees, we get the other ``half'' of
statement (a).\\

(b,c) Conversely, let $h$ be a minimal $P_{\al}$-harmonic function on
$\DL(q,r)$. By Theorem \ref{split-theorem}, $h(x_1x_2)=h_1(x_1) + h_2(x_2)$
with $h_i$ non-negative and $P_i$-harmonic. One of the $h_i\,$, say $h_1\,$, must
be positive. Minimality yields $h_1(x_1) = c_1\cdot h(x_1x_2)$ for
all $x_1x_2 \in \DL(q,r)$, where $c_1 > 0$. 
Thus $h(x_1x_2)$ depends only on $x_1$. 
Without loss of generality, $c_1 = 1$, and $h(x_1x_2) = h_1(x_1)$
for all $x_1x_2$. Minimality of $h$ with respect to $P_{\al}$ yields
minimality of $h_1$ with respect to $P_1$. Thus, by Proposition 
\ref{tree-martin}, $h(x_1x_2) = K_1(x_1,\xi_1)$ for some $\xi_1 \in \bd\T_q$.

The case when $h(x_1x_2)$ depends only on $x_2$ is analogous.\\

To complete the proof of statments (b) and (c), we have to study
minimality of the functions $x_1x_2 \mapsto K_i(\cdot,\om_i)$, $i=1,2$ 
with respect to $P_{\al}$.\\

If $\al = 1/2$ then by (\ref{martin-omega}), $K_i(\cdot,\om_i) \equiv 1$.
In this case it is known from {\sc Kaimanovich and Woess} \cite{KaiWoe}, 
\S 6.2, that the Poisson boundary is trivial, i.e., all bounded harmonic
functions are constant, which is the same as minimality of the constant function
$\mathbf 1$. This proves (c). \\

Suppose $\al \ne 1/2$. Then -- again by \cite{KaiWoe} --  the Poisson 
boundary is nontrivial, and the constant function $\mathbf 1$ is non-minimal.
If $\al > 1/2$, then this yields that $x_1x_2 \mapsto K_2(x_2,\om_2) = 1$
is non-minimal. On the other hand, we know from (\ref{martin-omega}) that
$$
g(x_1x_2):=K_1(x_1,\om_1) = 
\left(\frac{1-\al}{\al}\right)^{\mbox{\footnotesize $\hor(x_1)$}}.
$$
We can conjugate $P_{\al}$ by $g$, that is, we set
$$
\check p(x_1x_2,y_1y_2) = \frac{p(x_1x_2,y_1y_2)g(y_1y_2)}{g(x_1x_2)}\,.
$$
Then $g$ is minimal $P_{\al}$-harmonic if and only if $\mathbf 1$ is
minimal $\check P_{\al}$-harmonic. However, $\check P_{\al} = P_{1-\al}$
by a straightforward computation, and the constant function $\bf 1$ is not
minimal $P_{1-\al}$-harmonic by non-triviality of the Poisson boundary.
Thus, also $g$ is non-minimal for $P_{\al}$.
Again, the case $\al < 1/2$ follows by exchanging the roles of the
two trees.
\end{proof} 

We remark that for our nearest neighbour case, minimality, resp. 
non-minimalty of $\mathbf 1$ can be proved in a more elementary
(somewhat longer) way than by appealing to the results of 
{\sc Kaimanovich and Woess} \cite{KaiWoe}.\\

\begin{exa}\label{SRW-example}
We conclude this section by retranslating the results for simple random
walk on $\DL(q,q)$ to the setting and notation of the random walk
(\ref{walk-switch}) on $\DL(q,q)$ (``walk forward and switch or switch
and walk backward''). We write $\T^1$ and $\T^2$ for the first and  the
second tree, respectively. (Both are copies of $\T_q$.) We have $\al = 1/2$,
and the constant harmonic function $\mathbf 1$ is minimal.

In terms of configurations, each $\xi_1 \in \bd^*\T^1$ corresponds to an
infinite configuration
$$
\xi_1: \Z \to \Z_q \quad\text{with}\quad \xi_1(n) = 0 \; 
\forall\ n \le n_0(\xi_1) \in \Z\,.
$$
If we label the edges of $\T^1$ by elements of $\Z_q$, 
as described in (\ref{tree-seq}), then $\xi_1(n)$ is the label of
the edge between the horocycles $H_{n-1}(\T^1)$ and $H_n(\T^1)$ 
on the infinite geodesic
$\geo{\om_1\,\xi_1}$. Now let $(\eta,k) \in \Z_q \wr \Z$, and 
consider the associated pair $x_1x_2 \in \DL(q,q)$.
From the computations in Example \ref{alpha-walk}, we know that
$$
K_1(x_1,\xi_1) = 
q^{\mbox{\footnotesize $\hor(x_1 \cf \xi_1) - \hor(o_1 \cf \xi_1)$}}
$$
It is easy to compute in terms of $(\eta,k)$
\begin{equation}\label{positive-defect}
\begin{aligned}
\hor(x_1\cf\xi_1) &= \begin{cases} \min \{ n \le k : \xi_1(n+1) \ne \eta(n+1) \}
                                        & \text{if such $n$ exists,} \\
                                      k &\text{otherwise,}
                                \end{cases}
\\[4pt]
\hor(o_1\cf\xi_1) &= \begin{cases} \min \{ m \le 0 : \xi_1(m+1) \ne 0 \}
                                        & \text{if such $m$ exists,} \\
                                      0 &\text{otherwise.}
                                \end{cases}
\end{aligned}
\end{equation}
We shall write $\df^+\bigl((\eta,k),\xi_1) = \hor(x_1 \cf \xi_1) - \hor(o_1 \cf \xi_1)$,
the \emph{(positive) defect} of $(\eta,k)$ with respect to $\xi_1$.

Analogously, each $\xi_2 \in \bd^*\T^2$ corresponds to an
infinite configuration
$$
\xi_2: \Z \to \Z_q \quad\text{with}\quad \xi_2(n) = 0 \; 
\forall\ n \ge n_0(\xi_2) \in \Z\,.
$$
Here, we label the edges of $\T^2$ by elements of $\Z_q$, 
and $\xi_2(n)$ is the label of the edge between the horocycles
$H_{-n}(\T^2)$ and 
$H_{-n+1}(\T^2)$ on the infinite geodesic $\geo{\om_2\,\xi_2}$. 
If $(\eta,k) \in \Z_q \wr \Z \leftrightarrow x_1x_2 \in \DL(q,q)$, then 
$$
K_1(x_2,\xi_2) = 
q^{\mbox{\footnotesize $\hor(x_2 \cf \xi_2) - \hor(o_2 \cf \xi_2)$}}\,,
$$
and we compute
\begin{equation}\label{negative-defect}
\begin{aligned}
-\hor(x_2\cf\xi_2) &= \begin{cases} \max \{ n > k : \xi_2(n) \ne \eta(n) \}
                                       & \text{if such $n$ exists,} \\
                                     k &\text{otherwise,}
                                \end{cases}
\\[4pt]
-\hor(o_2\cf\xi_2) &= \begin{cases} \max \{ m > 0 : \xi_2(m) \ne 0 \}
                                       & \text{if such $m$ exists,} \\
                                     0 &\text{otherwise.}
                                \end{cases}
\end{aligned}
\end{equation}
We shall write $\df^-\bigl((\eta,k),\xi_2) = \hor(x_2 \cf \xi_2) - \hor(o_2 \cf \xi_2)$,
the \emph{(negative) defect} of $(\eta,k)$ with respect to $\xi_2$.

\vspace{.3cm}

We conclude:
write $\bd^+ (\Z_q \wr \Z)$ and $\bd^-(\Z_q \wr \Z)$ 
for all infinite configurations $\xi_1$, resp. $\xi_2$ as above.
(Every \emph{finitely supported} configuration appears once in each of the two parts of the
boundaries~!)

Then all  non-constant minimal $P$-harmonic functions are given by
\begin{equation}
\begin{aligned}
K_1(\cdot, \xi_1) &= 
q^{\mbox{\footnotesize $\df^+(\cdot,\xi_1)$}}\,,\;\xi_1 \in \bd^+(\Z_q\wr\Z)\,, 
\AND\\ K_2(\cdot,\xi_2) &= 
q^{\mbox{\footnotesize $\df^-(\cdot,\xi_2)$}}\,,\;\xi_2 \in \bd^-(\Z_q\wr\Z)\,.
\end{aligned}
\end{equation}
The constant function $\mathbf 1$ is also minimal harmonic.
\end{exa}

\medskip

\section{Switch-walk-switch}\label{extension}

We now turn our attention to the random walk (\ref{switch-walk-switch}), where at each step, the
lamplighter first switches the lamp at his actual position to a random state, then walks, and
then switches the lamp at the arrival point to a random state. 

As explained in Section \ref{geometry},
this corresponds to simple random walk on the modification of $\DL(q,q)$ where in the first
tree, additional edges are drawn between every vertex and the siblings of its predecessor, while
the second tree remains as it is.  

More generally, consider $\DL(q,r)$. For $x_1,u_1 \in \T_q$, we introduce the 
\emph{sibling relation} $x_1 \sib u_1 \!\!\iff\!\! x_1^-=u_1^-$.
We extend this relation to $\DL(q,r)$ by setting $x_1x_2 \sib u_1x_2 \!\!\iff\!\! x_1 \sib u_1$.
The new ege set on the same vertex set $\{ x_1x_2 \in \T_q \times \T_r : \hor(x_1)+\hor(x_2)=0\}$
is now given by
$$
\{ [x_1x_2, y_1y_2] : y_1^- \sib x_1 \;\text{and}\; y_2=x_2^- \}\,.
$$
We write $\DL^s(q,r)$ for the resulting graph.
Every vertex $x_1x_2$ with $\hor(x_1) = k$ has $q^2$ neighbours $y_1y_2$ 
with $\hor(y_1)=k+1$ and $qr$ neighbours with $\hor(y_1)=k-1$. 
Adapted to this structure, we choose $0 < \al < 1$ and consider the random walk 
on $\DL^s(q,r)$ with transition matrix $Q =Q_{\al}$ given by 
\begin{equation}\label{sws-walk}
q(x_1x_2,y_1y_2) = \begin{cases} \al/q^2 
                   & \text{if}\; y_1^- \sib x_1 \;\text{and}\;y_2=x_2^-\\
                                 (1-\al)/(qr) 
		   & \text{if}\; y_1 \sib x_1^- \;\text{and}\;y_2^-=x_2\\
                 0 & \text{otherwise,}
                   \end{cases}
\end{equation}
Now note that when $x_1 \sib u_1$ then transitions from $x_1x_2$ and $u_1x_2$ 
go to the same neighbours with the same probabilities. Thus, 
$Qh(x_1x_2) = Qh(u_1x_2)$ whenever $x_1 \sib u_1$, and we have
the following.

\begin{lem}\label{sibling-harmonic}
Every $Q_{\al}$-harmonic function is constant on each equivalence class of the sibling
relation.
\end{lem}

We can construct the \emph{factor graph} of $\DL^s(q,r)$ with respect to the sibling relation.
We write $[x_1]x_2$ for the equivalence class of $x_1x_2$, since all its elements
have the same second ``coordinate'', and $[x_1]$ is the sibling class in the first tree. 
Then the vertex set of the factor graph is $\{ [x_1]x_2 : x_1x_2 \in \DL^s(q,r)  \}$, and 
two classes $[x_1]x_2$ and $[y_1]y_2$ are connected by an edge of the factor graph
if and only if there is an edge bewteen a pair of representatives in $\DL^s(q,r)$.
Thus, if $\hor(y_2) = \hor(x_2) - 1$, there is an edge from $[x_1]x_2$ to $[y_1]y_2$
precisely when  $[y_1^-] = [x_1]$ and $x_2^- = y_2$. 
We write $\pi$ for the natural projection.
The next lemma is now immediate.

\begin{lem}\label{factor-graph}
 The factor graph of $\DL^s(q,r)$ with respect to the sibling relation is
$\DL(q,r)$. 

\vspace{.2cm}

The transition matrix $Q_{\al}$ is compatible with the factorization,
and its image under the projection $\pi: \DL^s(q,r) \to \DL(q,r)$
is $P_{\al}$, as defined in (\ref{random-walk}).
\end{lem}

By ``compatible'' we mean that 
$
q_{\al}(v_1x_2,[y_1]y_2) = \sum_{w_1 \sib y_1} q_{\al}(x_1x_2,w_1y_2)
$
is the same for each representative $v_1 \in [x_1]$, and this common value is the transition
probability from $[x_1]x_2$ to $[y_1]y_2$ of the projection of $Q_{\al}$.

\begin{cor}\label{factor-harmonic}
Every $Q_{\al}$-harmonic function is of the form $h \circ \pi$,
where $h$ is a $P_{\al}$-harmonic function on $\DL(q,r)$ and 
$\pi: \DL^s(q,r) \to \DL(q,r)$ is the factor map with respect to the 
sibling relation.
\end{cor}

\begin{exa}\label{sws-example}
Our final task is to retranslate once more to the lamplighter group,
by giving a direct description of the minimal harmonic functions for the
switch-walk-switch model that does not involve the above factor map.
We have $Q = Q_{1/2}$ on $\DL^s(q,q)$.

Now, it is clear what the factor map does to a pair $(\eta,k) \in \Z_q \wr \Z\,$:
it ``forgets'' (cancels) the value $\eta(k)$, and what remains is the
pair $(\eta_{\not \,k}, k)$, where
$$
\eta_{\not \,k}(n) = \begin{cases} \eta(n-1) & \text{if}\; n \le k\,,\\
                                 \eta(n)  & \text{if}\;n > k\,. 
                   \end{cases}
$$
Thus, with respect to the computations of Example \ref{SRW-example},
(\ref{negative-defect}) remains unchanged, while instead of the positive defect
we need
 \begin{equation}\label{new-positive-defect}
\begin{aligned}
\df^{\oplus}\bigl((\eta,k),\xi_1) &= \hor([x_1] \cf \xi_1) - \hor([o_1] \cf \xi_1)\,,
           \quad\text{where}\\[6pt]
\hor([x_1]\cf\xi_1) &= \begin{cases} \min \{ n \le k : \xi_1(n) \ne \eta(n) \}
                                        & \text{if such $n$ exists,} \\
                                      k &\text{otherwise,}
                                \end{cases}
\\[4pt]
\hor([o_1]\cf\xi_1) &= \begin{cases} \min \{ m \le 0 : \xi_1(m) \ne 0 \}
                                        & \text{if such $m$ exists,} \\
                                      0 &\text{otherwise.}
                                \end{cases}
\end{aligned}
\end{equation}
Again, the constant function $\mathbf 1$ is minimal $Q$-harmonic. 
All  non-constant minimal $Q$-harmonic functions are given by
\begin{equation}
\begin{aligned}
K_1(\cdot, \xi_1) &= 
q^{\mbox{\footnotesize $\df^{\oplus}(\cdot,\xi_1)$}}\,,\;\xi_1 \in \bd^+(\Z_q\wr\Z)\,, 
\AND\\ K_2(\cdot,\xi_2) &= 
q^{\mbox{\footnotesize $\df^-(\cdot,\xi_2)$}}\,,\;\xi_2 \in \bd^-(\Z_q\wr\Z)\,.
\end{aligned}
\end{equation}
\end{exa}

\medskip

\section{Final observations and speculations}\label{final}

\begin{imp} Even when $r \ne q$, one can interpret $\DL(q,r)$ as a ``lamplighter graph''
over $\Z\,$: at each point of $\Z$, there are green lamps with $q$ different states,
including ``off'', and red lamps with $r$ different states, again including ``off''.
The lamplighter walking along $\Z$ has to make sure that when his actual position
is $k \in \Z$, then the lamps in $(-\infty\,,\,k]$ have to be in one of the green
states, and those in $[k+1\,,\,+\infty)$ in one of the red states.
\end{imp}

\begin{imp} Several basic properties of random walks on $\DL(q,r)$ that are not
necessarily of nearest neighbour type, but invariant under a transitive group of
automorphisms of $\DL(q,r)$, were studied by {\sc Bertacchi} \cite{Ber}.
For a large class of random walks of this type, the Poisson boundary was 
determined by {\sc Kaimanovich and Woess} \cite{KaiWoe},
as mentioned above.
\end{imp}

\begin{imp} As a general principle, the three problems of (i) determining the Poisson
boundary ($\equiv$ all bounded harmonic functions), (ii) determining all minimal
harmonic functions  ($\equiv$ all positive harmonic functions), and (iii) determining
the full Martin compactification ($\equiv$ finding the directions of convergence
of the Martin kernels) have an increasing degree of difficulty. (As a matter of fact,
these three problems get sometimes mixed up even by advanced non-experts.)
Thus, the reader should not be a priori astonished by the fact that in this paper, we 
were able to solve problem (ii) for a much smaller class of random walks than those
for which (i) was solved in \cite{KaiWoe}.
\end{imp}

\begin{imp} In particular, to the author's knowledge, the present 
results provide the first example of a complete solution of problem (ii) 
on a finitely generated, solvable group. In addition, our group is
non-polycyclic.

On the other hand, the situation is much better understood for connected solvable
Lie groups, because more structure theory is at hand. For the basic example, namely
the \emph{affine group} over $\R$, random walks and harmonic functions were
studied in much detail, see {\sc Molchanov} \cite{Mol}, {\sc Elie} \cite{Eli},
or {\sc Bougerol and Elie} \cite{BouEli}. We recall here that the main result
of \cite{BouEli} implies existence of non-constant positive harmonic functions
for finite range random walks on finitely generated 
polycyclic groups with exponential growth, but a complete solution of (ii)
is not available for those groups.

As for the affine group over $\R$, the Poisson boundary
and the Martin compactification are well understood for random walks on the 
affine group over the $p$-adic numbers, see {\sc Cartwright, Kaimanovich and 
Woess} \cite{CarKaiWoe} and {\sc Brofferio} \cite{Bro}.
\end{imp}

\begin{imp} As pointed out by {\sc Bertacchi} \cite{Ber}, there is a natural geometric
compactification of $\DL(q,r)$. Namely, this graph is a subgraph of the direct
product $\T_q  \times \T_r$, for which a natural compactification is 
$\wh \T_q \times \wh \T_r$. Thus, we define $\wh\DL(q,r)$ as the closure of
$\DL(q,r)$ in $\wh\T_q \times \wh\T_r$, and 
$\bd\DL(q,r) = \wh\DL(q,r) \setminus \DL(q,r)$. Almost sure convergence 
of random walks to a boundary-valued random variable is studied in \cite{Ber}.
However, at present, we are far from proving that this compactification is
in some sense comparable or identical with the Martin compactification even in the
case of simple random walk.
\end{imp}

\begin{imp} The following seems noteworthy regarding the two classes of examples that
we have studied: in the case of drift ($\al \ne 1/2$), the minimal harmonic functions
are parametrized (continuously on each part) by $\bd^*\T_q \cup \bd^*\T_r$.
In the drift-free case ($\al = 1/2$), the \emph{additional} minimal harmonic function
$\mathbf 1$ enters the stage. Thus, in some sense, the cone of positive harmonic
functions is \emph{bigger} in the driftfree case than in presence of drift. This 
contrasts with all examples known so far. (Of course, the constant 
function $\mathbf 1$ is a Martin kernel even when $\al \ne 1/2$, but then 
it does not belong to $\MM_{\min}$.)
\end{imp}

\begin{imp} Finally, the reason why our method does not adapt to the ``walk or switch''
model (\ref{walk-or-switch}) appears to be that in contrast with the cases that we have
solved here, this random walk is \emph{not} invariant under a group of automorphisms
of $\DL(q,r)$ that acts \emph{transitively} both on $\bd^*\T^1$ and $\bd^*\T^2$.

A next  step could be to try to prove the Decomposition Theorem \ref{split-theorem} for all
irreducible random walks with the latter transitivity property.
\end{imp}

\noindent{\bf Acknowledgement.} The author is grateful to R\"ognvaldur I. M\"oller for
pointing out the relationship between the Diestel-Leader graphs $\DL(q,q)$ and 
the lamplighter groups.

\end{document}